\documentclass[12pt]{amsart}

\usepackage{epsfig}

\usepackage{amsmath}
\usepackage{amssymb,amsthm}
\usepackage{graphicx}
\usepackage{dsfont}

\usepackage{hyperref}

\usepackage{enumerate}

\usepackage{pst-node}

\usepackage{tikz}

\usepackage{a4wide}

\newtheorem{propo}{Proposition}[section]

\newtheorem{lemma}[propo]{Lemma}

\newtheorem{theo}[propo]{Theorem}

\newtheorem{cor}[propo]{Corollary}
\newtheorem{prop}[propo]{Proposition}
\newtheorem{rem}[propo]{Remark}

\def\Cay{{\rm Cay}}

\def\qed{\ifmmode\square\else\nolinebreak\hfill
$\square$\fi\par\vskip12pt}

\usepackage{indentfirst,latexsym,bm}

\renewcommand{\proof}{\par\noindent{\bf Proof.\ \ }}

\def\Sym{{\rm Sym}}

\def\Aut{{\rm Aut}}

\def\diam{{\rm diam}}

\def\K{{\rm K}}

\begin{document}
\title{Finite $s$-geodesic-transitive digraphs}

%\thanks{Supported by the NSF of China (11561027,11661039,61563018) and NSF of Jiangxi (GJJ150460, GJJ150444,20151BAB201001).}

\thanks{Supported by the NNSF of China (12271524,12061034),   NSF of Jiangxi (20224ACB201002, 20212BAB201010) and NSF of Hunan (2022JJ30674)}

\author[W. Jin]{Wei Jin}
\address{Wei Jin}
\address{School of Statistics, Key Laboratory of Data Science in Finance and Economics\\
Jiangxi University of Finance and Economics\\
 Nanchang, Jiangxi, 330013, P.R.China}
 \address{School of Mathematics and Statistics\\
Central South University\\
Changsha, Hunan, 410075, P.R.China}
\email{jinweipei82@163.com}
%\email{jinwei@jxufe.edu.cn}

%School of Mathematics and Statistics, HNP-LAMA, Central South University,
%  Changsha, Hunan 410083, P. R. China

%\date{ }
%\date\today

\maketitle

\begin{abstract}
This paper  initiates the investigation of  the family of  $(G,s)$-geodesic-transitive
digraphs with  $s\geq 2$. We first  give a global analysis   by  providing a reduction result. Let $\Gamma$ be  such a digraph and let   $N$ be a  normal subgroup of $G$ maximal with respect to having at least $3$ orbits.   Then the quotient digraph $\Gamma_N$ is
$(G/N,s')$-geodesic-transitive where $s'=\min\{s,\diam(\Gamma_N)\}$,  $G/N$ is  either quasiprimitive or bi-quasiprimitive on $V(\Gamma_N)$, and $\Gamma_N$ is either directed or an  undirected complete graph.
Moreover, it is further shown that if
$\Gamma $ is not $(G,2)$-arc-transitive,  then $G/N$ is   quasiprimitive  on $V(\Gamma_N)$.

On the other hand, we also consider the  case that the normal subgroup $N$ of $G$ has one orbit on the vertex set. We show that if $N$ is regular  on $V(\Gamma)$, then $\Gamma$ is a circuit, and particularly each  $(G,s)$-geodesic-transitive
normal Cayley digraph with   $s\geq 2$,   is a circuit.

Finally, we investigate  $(G,2)$-geodesic-transitive digraphs with either valency at most 5 or diameter at most 2.
Let  $\Gamma$ be a $(G,2)$-geodesic-transitive
digraph. It is proved that: if $\Gamma$ has   valency at most $5$, then  $\Gamma$ is  $(G,2)$-arc-transitive; if   $\Gamma$ has  diameter $2$, then $\Gamma$  is  a balanced incomplete block design with the Hadamard parameters.

\end{abstract}

\vspace{2mm}
\vspace{2mm}

 \hspace{-17pt}{\bf Keywords:}   $s$-geodesic-transitive digraph, automorphism group, permutation group.

 \hspace{-17pt}{\bf Math. Subj. Class.:} 05E18; 20B25

\section{Introduction}

A finite \emph{digraph} (short for \emph{directed graph}) $\Gamma$ consists of a finite set $V(\Gamma)$ of vertices and an antisymmetric irreflexive
relation $\rightarrow$ on $  V(\Gamma)$.  An \emph{arc}  of $\Gamma$ is an ordered pair of adjacent vertices. Hence for two vertices $u$ and $v$,
$u\rightarrow v$ is equivalent to that $(u,v)$ is an arc.
For each vertex  $v\in  V(\Gamma)$, we use  $\Gamma^-(v) = \{u\in V|u\rightarrow v\}$ to denote the set of \emph{in-neighbours} of $v$ and use
$\Gamma^+(v) = \{u\in V|v\rightarrow u\}$ to denote the set of \emph{out-neighbours} of $v$.
We say that a digraph $\Gamma$ is $k$-regular if both the set $\Gamma^-(v) $  and the set $\Gamma^+(v) $  have
size $k$ for all $v\in  V(\Gamma)$, and  $\Gamma$ is \emph{regular} if it is $k$-regular for some positive
integer $k$.

For a non-negative integer $s$, an $s$-arc
in a digraph (or graph) $\Gamma$ is a sequence $(v_0, v_1, \ldots, v_s)$ of vertices with $v_i\rightarrow v_{i+1}$ for each $i = 0,\ldots, s - 1$.
This digraph (or graph)    is said to be  \emph{ $(G,s)$-arc-transitive} if its automorphism subgroup $G$  is transitive on
all the $s$-arcs.
The family of $s$-arc-transitive undirected graphs has been studied intensively, beginning with the seminal result of Tutte \cite{Tutte-1,Tutte-2},
and particularly  it is shown in 1981 by Weiss \cite{weiss} that finite undirected graphs of valency at least 3 can only be $s$-arc-transitive for $s\leq 7$,
more work  refer to \cite{ACMX-1996,GLP-1,IP-1,MW-2000,Praeger-4,Weiss-1}.
In  contrast with the situation for undirected graphs,  there exist  infinitely many classes of $s$-arc-transitive digraphs with unbounded $s$ other than directed cycles. Constructions for such classes of digraphs were initiated by Praeger  \cite{Praeger-1989} in 1989  and have stimulated a lot of research.
A few years later, Conder, Lorimer and Praeger
\cite{CLP-1995}   showed that  for
every integer $k \geq 2$ and every integer $s \geq 1$ there are infinitely many finite $k$-regular $(G, s)$-arc-transitive digraphs with $G$ quasiprimitive on the vertex set.
Giudici, Li and Xia \cite{GLX-2017} in 2017 solved the long-standing existence problem of vertex-primitive 2-arc-transitive digraphs by constructing an infinite class of such digraphs.
After one year, Giudici and Xia \cite{GX-2018} investigated  vertex-quasiprimitive 2-arc-transitive digraphs, and reduced the problem of
vertex-primitive 2-arc-transitive digraphs to almost simple groups, and it  includes a complete
classification of vertex-quasiprimitive 2-arc-transitive digraphs where the action on
vertices has quasiprimitive type SD or CD.
Given integers $k$ and $m$, Morgan,  Poto\v cnik and  Verret \cite{MPV-2019}
 constructed a $G$-arc-transitive graph of valency $k$ and an $L$-arc-transitive oriented digraph of out-valency $k$ such that $G$ and $L$ both admit blocks of imprimitivity of size $m$, and
more work see \cite{CPW-1993,Evans-1997,MMSZ-2002,WS-2003}.

%A  2-arc $(u,v,w)$ is called a \emph{$2$-geodesic} if it satisfying that the distance from $u$ to $w$ is 2, that is, $u\rightarrow v \rightarrow w$ but  $u \nrightarrow w$.
%A 2-geodesic $(u,v,w)$ is a \emph{T-$2$-geodesic} if further  $w \rightarrow u$; and a 2-geodesic $(u,v,w)$ is a \emph{NT-$2$-geodesic} if   $w \nrightarrow u$.
%Obviously, if a digraph has both T-$2$-geodesics and NT-$2$-geodesics, then it is not 2-geodesic-transitive.

An $s$-arc $(v_0, v_1, \ldots, v_s)$ in a digraph $\Gamma$ is called an \emph{$s$-geodesic} if the distance from $v_0$ to $v_s$ is $s$.
A  digraph (or graph)    is said to be  \emph{ $(G,s)$-geodesic-transitive} if its automorphism subgroup $G$  is transitive on
 the set of  $i$-geodesics for each $i\leq s$.
The simplest of $s$-geodesic-transitive digraphs are  directed cycles.

By definition, each $s$-geodesic of  $\Gamma$ is an $s$-arc,
but the converse is not true, for instance a 2-arc $(u,v,w)$ of $\Gamma$ satisfying $u\rightarrow w$ is  not a 2-geodesic. Thus the family
of  $s$-arc-transitive digraphs is  contained in
the family of $s$-geodesic-transitive digraphs.

The possible local structures of  $s$-geodesic-transitive undirected graphs for $s\geq 2$ are
characterized by Devillers, Li, Praeger and the author \cite{DJLP-clique}, it is proved that
for a vertex $u$,  either   $[\Gamma(u)]\cong
m\K_r$ for some  integers $m\geq 2,r\geq 1$, or
 $[\Gamma(u)]$ is a connected graph of diameter $2$.
And the families of $2$-geodesic-transitive graphs of valency 4 and of prime valency have been
determined in \cite{DJLP-2} and \cite{DJLP-prime}, respectively.
The 2-geodesic-transitive
undirected graphs have also been extensively studied in the literature, see
for example,
\cite{FH-2018,HFZ-1,DJLP-compare}.

In this paper, we initiate the study of finite $(G,s)$-geodesic-transitive digraphs for $s\geq 2$.
Our first theorem gives a global analysis  for such  digraphs, and it provides a reduction result.

\begin{theo}\label{2gtdig-quotient-1}
Let $\Gamma $ be a connected $(G,s)$-geodesic-transitive
digraph for some  $s\geq 2$, and let   $N$ be a  normal subgroup of $G$ maximal with respect to having at least $3$ orbits.   Then $\Gamma_N$ is connected
$(G/N,s')$-geodesic-transitive where $s'=\min\{s,\diam(\Gamma_N)\}$, and $G/N$ is  either quasiprimitive or bi-quasiprimitive on $V(\Gamma_N)$.
Moreover, $\Gamma_N$ is either a digraph  or an  undirected complete graph.
\end{theo}

Theorem \ref{2gtdig-quotient-1} directly leads to the following corollary which shows that for  $(G,2)$-geodesic-transitive
but not $(G,2)$-arc-transitive digraphs, if $G$ has  a  normal subgroup $N$ maximal with respect to having at least $2$ orbits, then $G/N$ is not   bi-quasiprimitive on $V(\Gamma_N)$.

\begin{cor}\label{2gtdig-quotient-cor1}
Let $\Gamma $ be a connected $(G,2)$-geodesic-transitive
but not $(G,2)$-arc-transitive digraph. Let   $N$ be a  normal subgroup of $G$ maximal with respect to having at least $2$ orbits.   Then
$G/N$ is   quasiprimitive  on $V(\Gamma_N)$, and  $\Gamma_N$ is either a connected
$(G/N,2)$-geodesic-transitive digraph or a  connected
$G/N$-arc-transitive  undirected complete graph.
\end{cor}

By  Theorem  \ref{2gtdig-quotient-1},  for $s\geq 2$ each connected $(G, s)$-geodesic-transitive digraph has a connected  $(G, s)$-geodesic-transitive quotient digraph
corresponding to a normal subgroup $N$ of $G$ such that $G/N$ acts quasiprimitively or bi-quasiprimitively on the vertex set of the quotient digraph. Thus a preliminary step in determining all $s$-geodesic-transitive digraphs may be the determination of the base one.
We  achieved this in  Proposition \ref{2gtdig-irreducible-2},
in the case where $G$ is soluble. It is shown that
 $\Gamma $ is a circuit with $r$ vertices  where $r = 4$ or
$r$ is a prime.
Our next  theorem is another  contribution to the determination of
 connected $(G,s)$-geodesic-transitive  vertex quasiprimitive digraphs.

For integer   $s\geq 2$ and  a $(G,s)$-geodesic-transitive
digraph $\Gamma $, the  second   theorem shows that
if $G$ has a  nontrivial regular normal subgroup, then $\Gamma$ is known.

\begin{theo}\label{2gtdig-regularnormal-1}
Let $\Gamma $ be a $(G,s)$-geodesic-transitive
digraph for some  $s\geq 2$, and let  $N$ be a
nontrivial normal subgroup of $G$.  Suppose that  $N$ is regular  on $V(\Gamma)$.
Then $\Gamma$ is a circuit. In particular, each  $(G,s)$-geodesic-transitive
$G$-normal Cayley digraph with   $s\geq 2$,   is a circuit.

\end{theo}

We give a remark of Theorem \ref{2gtdig-regularnormal-1}.

\begin{rem}\label{regular-rem-1}
{\rm  Let $\Gamma $ be a connected $(G,2)$-geodesic-transitive
but not $(G,2)$-arc-transitive digraph. If  $G$ is    quasiprimitive  on $V(\Gamma_N)$ of type HA, HS, HC or  TW, then $G$ has a normal subgraph that acts regularly on the vertex set, and  by Theorem \ref{2gtdig-regularnormal-1},
$\Gamma$ is a circuit.
}
\end{rem}

Our third  theorem investigates $(G,2)$-geodesic-transitive digraphs with either valency at most 5 or diameter at most 2.
%The techniques for investigating 2-geodesic-transitive digraphs are  quite different from those for
%$2$-geodesic-transitive (undirected) graphs.

\begin{theo}\label{2gtdig-smallval-th1}
Let  $\Gamma$ be a $G$-arc-transitive
digraph.  Then the following statements hold.

 \begin{itemize}
\item[(i)] If $\Gamma$ has   valency at most $5$, then $\Gamma$ is  $(G,2)$-geodesic-transitive if and only if $\Gamma$ is  $(G,2)$-arc-transitive.

\item[(ii)] If   $\Gamma$ is $(G,2)$-geodesic-transitive  of diameter $2$, then  $\Gamma$  is  a balanced incomplete block design with the Hadamard parameters.

\end{itemize}

\end{theo}

\section{Preliminaries}

In this section, we will give some definitions about groups and digraphs that will be used in the
paper. For the group theoretic terminology not defined here we refer the reader to \cite{Cameron-1,DM-1,Wielandt-book}.

All digraphs in this paper are finite and  simple. For a digraph $\Gamma$, we use $V(\Gamma)$ and $Arc(\Gamma)$ to denote
its vertex set and arc set, respectively.

An \emph{automorphism} of a digraph $\Gamma$ is a permutation $\pi$ of $V(\Gamma)$ which has the
property that $u\rightarrow v $ if and only if $u^\pi \rightarrow v^\pi$.
The set of all automorphisms
of $\Gamma$, with the operation of composition forms a group which   is
called the \emph{automorphism group} of $\Gamma$, and denoted by $\Aut(\Gamma)$.

A digraph $\Gamma$ is called \emph{$G$-vertex-transitive} or \emph{$G$-arc-transitive} if its automorphism subgroup $G$ acts transitively on its vertex set or arc set, respectively. It is obvious  that each $G$-arc-transitive digraph is $G$-vertex-transitive and each $G$-vertex-transitive  digraph is regular.

The number of arcs traversed in the shortest directed path from  $u$ to $v$ is called the \emph{distance} in $\Gamma$ from
$u$ to $v$, and is denoted by $d_{\Gamma}(u,v)$. The maximum value of the distance function in $\Gamma$ is called the \emph{diameter} of $\Gamma$,
 and denoted by $\diam(\Gamma)$. Define $\Gamma_i^+(u)=\{v\in V(\Gamma)|d_{\Gamma}(u,v)=i\}$ for $i\geq 1$.
In particular $\Gamma_1^+(u)=\Gamma^+(u)$.

A subdigraph $X$ of a digraph $\Gamma$ is an \emph{induced subdigraph} if
$(u,v)$ is an arc of  $X$   if and only if $(u,v)$ is an arc  in $\Gamma$.  When $U\subseteq V(\Gamma)$, we denote by
$[U]$  the
subdigraph of $\Gamma$ induced by $U$.
Let $\Sigma$ be a digraph. For a positive integer $m$,  the digraph consisting of $m$ vertex disjoint copies of $\Sigma$ is denoted by $m\Sigma$.

A connected digraph $\Gamma$ is called  \emph{strongly connected} if,
for all $u,v\in V(\Gamma)$, there is a $t$-arc $(u=u_0, u_1,\ldots,u_t=v)$  for some positive integer $t$.
Let $\Gamma$ be a finite $G$-arc-transitive digraph. By \cite[Lemma 2]{Neumann-1977}, if   the underlying undirected graph of $\Gamma$ (with $\{u,v\}$ an edge if either $(u,v)$ or
$(v,u)$  is an arc of $\Gamma$) is connected, then  $\Gamma$ is  strongly connected.

%Arc-transitivity in the digraph implies that if two of the edges are of the form
%$(u,v)$ and $(v,u)$, then all the edges come in pairs. Hence, in order to study
%something new, we will restrict ourselves to arc-transitive digraphs with the
%property that if $(u, v) \in  Arc(\Gamma)$, then $(v, u) \notin  Arc(\Gamma)$.

The \emph{girth} or \emph{directed girth} of a digraph is the minimum length
of a closed path with at least three vertices.

%The girth for directed graphs is usually defined as the directed girth, that is, the minimum length of a directed cycle (or $\infty$ if no directed cycles exist).

For integers $r\geq 3$, an $r$-arc  $(w_0,w_1,\ldots,w_{r})$ with distinct vertices is called a \emph{circuit} of length $r$ if $w_{r}= w_0$.
A shortest circuit is called a \emph{minimal circuit}. Thus the girth  of $\Gamma$ is the length of a minimal
circuit.

A transitive permutation group $G\leqslant \Sym(\Omega)$ is said to be \emph{regular} on $\Omega$, if for any $\omega\in \Omega$, the stabilizer $G_{\omega}=1$.

For a non-empty subset $S$ of a group $H$ the \emph{Cayley digraph}
$\Gamma = \Cay(H, S)$ is defined to be the digraph with vertex set $H$ and with arc set $Arc(\Gamma) =
\{(h, xh)|h\in  H, x \in  S\}$. The arc set $Arc(\Gamma)$ is anti-symmetric provided that $S\cap  S^{-1} $ is empty,
where $S^{-1} = \{x^{-1}|x\in S\}$, $\Gamma$ is regular of valency $|S|$, and $\Gamma$ is connected whenever $S$ is
a generating set for $H$. Also, $\Gamma$ admits as an automorphism group the semidirect product
$H:\Aut (H,S)$, where $H$ acts by right translation and the set stabilizer $\Aut (H,S)=\{\alpha\in \Aut(H)| S^\alpha=S\}$ of $S$ in the
automorphism group of $H$ acts by conjugation; in particular, $\Aut (H,S)$ acts transitively on
$S$ by conjugation if and only if $\Cay(H, S)$ is $(H: \Aut (H,S), 1)$-arc transitive. Moreover,
$\Cay(H, S)$ is said to be \emph{$G$-normal} if $H$ is a normal subgroup of  $G$.

Let $\Gamma$ be a connected $(G, s)$-geodesic-transitive digraph where $s\geq  1$, and let $N$ be
a normal subgroup of $G$ with more than two orbits in $V(\Gamma)$.
 Then  the \emph{quotient digraph} $\Gamma_N$ is defined as the digraph with
vertices the $N$-orbits in $V(\Gamma)$ and with $(A, B)$ an arc, where $A$ and $B$ are $N$-orbits, if and only
if, for some $a\in A$  and $b\in B$, $(a,b)$ is an arc of $\Gamma$.

A transitive permutation group $G\leqslant \Sym(\Omega)$ is said to be \emph{quasiprimitive}, if every non-trivial normal subgroup of $G$ is transitive on $\Omega$, while
$G$ is said to be \emph{bi-quasiprimitive} if every non-trivial normal subgroup of $G$ has at most two orbits on $\Omega$ and there exists one which has exactly two orbits on $\Omega$.
Quasiprimitivity  is a generalization
of primitivity as every normal subgroup of a primitive group is transitive,
but there exist quasiprimitive groups which are not primitive.  Praeger \cite{Praeger-4} generalized the O'Nan-Scott Theorem for primitive groups to quasiprimitive
groups and showed that a finite quasiprimitive group is one of eight distinct types: Holomorph Affine (HA), Almost Simple (AS), Twisted Wreath product (TW), Product Action (PA), Simple Diagonal (SD), Holomorph Simple (HS), Holomorph Compound (HC) and Compound Diagonal (CD). For more information about quasiprimitive and bi-quasiprimitive permutation groups, refer to \cite{Praeger-2003-biq}.

%Moreover, if $\Gamma$ is regular and
%$(G, s)$-arc-transitive with $s \geq 2$ then it is also $(G,s -1)$-arc-transitive.

%If $\Gamma$ is regular and
%$(G, s)$-geodesicc-transitive with $s \geq 2$,  is it also $(G,s -1)$-geodesic-transitive.

%\subsection{Some results}

\begin{lemma}\label{2gtdig-lem-at-1}
Let $\Gamma$ be a connected $G$-arc-transitive digraph of valency $k\geq 2$.   Let $(u,v)$ be an arc of $\Gamma$. Then the following statements hold.

 \begin{itemize}
\item[(1)]      $ \Gamma^+(u)\neq \{v\}\cup (\Gamma^+(u)\cap \Gamma^+(v))$.

\item[(2)]    $ \Gamma^+(u)\cap \Gamma^+(v)=\emptyset$
if and only if,   each $2$-arc of $\Gamma$ is a  $2$-geodesic.

\end{itemize}

\end{lemma}
\proof (1) Suppose that $ \Gamma^+(u)= \{v\}\cup (\Gamma^+(u)\cap \Gamma^+(v))$.
Then since $\Gamma$ is  $G$-arc-transitive,  the vertex stabilizer $G_u$ acts transitively on  $ \Gamma^+(u)$, and so  $[\Gamma^+(u)]$ is a digraph with out-valency $k-1$.
Therefore, for any two vertices $x,y \in \Gamma^+(u)$, we have two arcs $(x,y)$ and $(y,x)$.
Again using the    $G$-arc-transitive property of $\Gamma$,  we know that  $\Gamma$ must be  an undirected graph, which is a contradiction.
Thus $ \Gamma^+(u)\neq \{v\}\cup (\Gamma^+(u)\cap \Gamma^+(v))$.

(2) It is obvious. \qed

\section{Reduction}

Let $\Gamma$ be a connected $(G, s)$-geodesic-transitive digraph for some $s\geq 2$. In this section we study the nature of intransitive normal subgroups $N$ of $G$. The first lemma shows that each orbit of $N$ is arc-less.

\begin{lemma}\label{2gtdig-normalsubg-1}
Let $\Gamma $ be a connected $(G,s)$-geodesic-transitive
digraph for some  $s\geq 1$, and let  $N$ be a
nontrivial intransitive normal subgroup of $G$.
Then there is no $N$-orbit contains any arc of  $\Gamma$.

\end{lemma}
\proof Suppose that there exists one $N$-orbit $B_0$ that contains an arc $(u,v)$ of  $\Gamma$.
The  subgroup  $N$ is not transitive on $V(\Gamma)$ leading to that it  has at least two orbits. Since $\Gamma$
is connected and $N$ is transitive on $B_0$, it follows that $\Gamma^+(u)$ intersects nontrivially with some other $N$-orbit, say $B_1$, and set
$v'\in \Gamma^+(u)\cap B_1$. Then $(u,v')$ is an arc.
Since $\Gamma $ is  a $(G,s)$-geodesic-transitive
digraph for some  $s\geq 1$, $G_u$ has an element that can maps $v$ to $v'$, which is impossible, as $G_u$
fixes $B_0$ setwise. Thus there is no $N$-orbit contains any arc of  $\Gamma$.
\qed

\begin{lemma}\label{2gtdig-normalsubg-bip}
Let $\Gamma $ be a $(G,s)$-geodesic-transitive
digraph for some  $s\geq 2$, and let  $N$ be a
nontrivial  normal subgroup of $G$.
If  $N$ has $2$ orbits on $V(\Gamma)$, then  $\Gamma$ is $(G,2)$-arc-transitive and bipartite.

\end{lemma}
\proof
Suppose that $N$ has $2$ orbits on $V(\Gamma)$, say $\Delta_1$ and $\Delta_2$.
Then by Lemma \ref{2gtdig-normalsubg-1}, neither $\Delta_1$ nor $\Delta_2$ contains any arc of $\Gamma$, and so $\Gamma$ is a  bipartite digraph.

Let $(u,v)$ be an arc where $u\in \Delta_1$ and $v\in \Delta_2$.
Suppose that $\Gamma^+(u)\cap \Gamma^+(v)\neq \emptyset$. Let
$v'\in \Gamma^+(u)\cap \Gamma^+(v)$. Then $(u,v')$ and $(v,v')$ are two arcs. Moreover, $v'$ is either in $\Delta_1 $ or in $ \Delta_2$.
If $v'\in \Delta_1$, then $\Delta_1$ contains the arc $(u,v')$, a contradiction.
If $v'\in \Delta_2$, then $\Delta_2$ contains the arc $(v,v')$, again a contradiction.
Thus $\Gamma^+(u)\cap \Gamma^+(v)= \emptyset$, and so by Lemma \ref{2gtdig-lem-at-1} each 2-arc of $\Gamma$ is a 2-geodesic.
It concludes that   $\Gamma$ is $(G,2)$-arc-transitive. \qed

Let $\Gamma_N$ be the quotient digraph with
vertices the $N$-orbits in $V(\Gamma)$ and with $(A, B)$ an arc, where $A$ and $B$ are $N$-orbits, if and only
if, for some $a\in A$  and $b\in B$, $(a,b)$ is an arc of $\Gamma$.

Now we show that each normal subgroup $N$ of $G$ with at least 3  orbits on vertices corresponds to a connected
 $(G, s)$-geodesic-transitive quotient digraph $\Gamma_N$.

\begin{lemma}\label{2gtdig-quotient-2orbit}
Let $\Gamma $ be a connected $(G,s)$-geodesic-transitive
digraph for some  $s\geq 2$, and let   $N$ be a  normal subgroup of $G$ with at least $3$ orbits.   Then $\Gamma_N$ is either directed or an  undirected complete graph.

\end{lemma}
\proof  Suppose that $\Gamma_N$ is not directed. Then  $\Gamma_N$ has two arcs $(B_1,B_2)$ and $(B_2,B_1)$ where $B_1,B_2$ are $N$-orbits. By the definition of $\Gamma_N$,  there exist $u,u'\in B_1$ and $v,v'\in B_2$
such that  $(u,v)$ and $(v',u')$ are two arcs of $\Gamma$. Since $N$ is transitive on each orbit, we can assume that $u=u'$. Then $(v,u,v')$ is a 2-arc of $\Gamma$.
By Lemma \ref{2gtdig-normalsubg-1}, neither $B_1$ nor $B_2$ contains an arc of $\Gamma$. Thus $(v,v')$ and $(v',v)$ are not arcs and hence  $(v,u,v')$ is a 2-geodesic  of $\Gamma$.

Since $\Gamma$ is  $(G,2)$-geodesic-transitive, it follows that $G_{v,u}$ is transitive on $\Gamma_2^+(v)\cap \Gamma^+(u)$.
The stabilizer  $G_{v,u}$ setwise fixes $B_1$ and $B_2$,   and so $\Gamma_2^+(v)\cap \Gamma^+(u)\subseteq B_2$.

Let $(B_1,B_3)$ be an arc of $\Gamma_N$ where $B_3(\neq B_2)$ is an $N$-orbit. Then   $(B_2,B_1,B_3)$ is a 2-arc of $\Gamma_N$, and
there exist $u''\in B_1$ and $w\in B_3$ such that $(u'',w)$ is an arc of $\Gamma$. The group   $N$ is transitive on each orbit leading to that  we can assume  $u=u''$. Hence  $(v,u,w)$ is a 2-arc.
Due to $\Gamma_2^+(v)\cap \Gamma^+(u)\subseteq B_2$ and $B_3\neq B_2$, we know  that $(v,w)$ must be an arc of $\Gamma$, and so
$(B_2,B_3)$ is an arc of $\Gamma_N$.
As a consequence  $\Gamma_N^+(B_1)=\{B_2\}\cup (\Gamma_N^+(B_1)\cap \Gamma_N^+(B_2))$
and $\Gamma_N^+(B_2)=\{B_1\}\cup (\Gamma_N^+(B_1)\cap \Gamma_N^+(B_2))$, and so  $(\Gamma_N)_2^+(B_1)\cap \Gamma_N^+(B_2)=\emptyset$.
Since $\Gamma_N$ is $G/N$-arc-transitive, it follows that $\Gamma_N$ has diameter 1 and it is an undirected complete  graph.
\qed

We are ready to prove Theorem \ref{2gtdig-quotient-1}.

\medskip
\noindent {\bf Proof of Theorem \ref{2gtdig-quotient-1}.}
Since $\Gamma $ is a connected digraph, it is easy to see that the quotient digraph $\Gamma_N $ is also connected.
Since $N$ is a  normal subgroup of $G$ maximal with respect to having at least $3$ orbits, it follows that all normal subgroups of $G/N$ are transitive or have two orbits on $V(\Gamma_N)$. Thus $G/N$ is quasiprimitive or bi-quasiprimitive on $V(\Gamma_N)$. Moreover, by Lemma \ref{2gtdig-quotient-2orbit}, $\Gamma_N$ is either directed or an  undirected complete graph.

 Let  $(B_0,B_1,B_2,\ldots,B_{t})$ and
$(C_0,C_1,C_2,\ldots,C_{t})$ be two $t$-geodesics of $\Gamma_N$
where $t\leq s'=\min\{s,\diam(\Gamma_N)\}$.

Then by the definition of  $\Gamma_N$, there exist $x_i\in B_i$ and $x_{i+1}'\in B_{i+1}$ such
that $(x_i,x_{i+1}')$ is an arc of $\Gamma$.
Since $N$ is transitive on each $B_i$, we have $x_j\in B_j$ such that
$(x_0,x_1,x_2,\ldots,x_{t})$ is a  $t$-geodesic of $\Gamma$. Similarly,
there exist $y_i\in C_i$ such that $(y_0,y_1,y_2,\ldots,y_{t})$
is a  $t$-geodesic of $\Gamma$.
As $t\leq s'\leq s$ and $\Gamma$
is $(G,s)$-geodesic-transitive,  the group $G$ has an element $g$ such that
$(x_0,x_1,x_2,\ldots,x_{t})^g=(y_0,y_1,y_2,\ldots,y_{t})$, and hence
$g$ induces an element $g'$ of $G/N$, such that $(B_0,B_1,B_2,\ldots,B_{t})^{g'}=(C_0,C_1,C_2,\ldots,C_{t})$. Thus
$\Gamma_N$ is $(G/N,s')$-geodesic-transitive. We conclude the proof. \qed

Theorem \ref{2gtdig-quotient-1} directly leads to Corollary \ref{2gtdig-quotient-cor1} which is a reduction result on $(G,2)$-geodesic-transitive
but not $(G,2)$-arc-transitive digraphs.

\medskip
\noindent {\bf Proof of Corollary \ref{2gtdig-quotient-cor1}.}
If   $N$ has exactly $2$ orbits on $V(\Gamma)$,
then by Lemma \ref{2gtdig-normalsubg-bip},   $\Gamma$ is $(G,2)$-geodesic-transitive  leading to that  $\Gamma$ is $(G,2)$-arc-transitive,  contradicts to our assumption.
Thus   $N$ has at least $3$ orbits on $V(\Gamma)$.
It follows from  Theorem \ref{2gtdig-quotient-1} that  $G/N$ is quasiprimitive or bi-quasiprimitive on $V(\Gamma_N)$, and  $\Gamma_N $ is connected  $(G/N,s')$-geodesic-transitive where $s'=\min\{2,\diam(\Gamma_N)\}$.
Moreover, as each arc-transitive digraph has diameter at least 2, we know that  $\Gamma_N$ is either a connected
$(G/N,2)$-geodesic-transitive digraph or a  connected
$G/N$-arc-transitive  undirected complete graph

Assume that $G/N$ acts  bi-quasiprimitively on $V(\Gamma_N)$.
Then $G/N$ has a nontrivial normal subgroup     which  has exactly $2$ orbits on $V(\Gamma_N)$, say $\Delta_1'$ and $\Delta_2'$.
Applying  Lemma \ref{2gtdig-normalsubg-1}, inside of  $\Delta_1'$ and  $\Delta_2'$ do not have any arc of $\Gamma_N$ and $\Gamma_N$ is a bipartite digraph.
It follows that
neither $\Delta_1'$ nor $\Delta_2'$ contains an arc of $\Gamma$, and so
 $\Gamma$ is a bipartite digraph. It concludes that  each 2-arc of $\Gamma$ is a 2-geodesic.
Therefore, $\Gamma$ is $(G,2)$-arc-transitive, a contradiction.
Hence  $G/N$ is not bi-quasiprimitive on $V(\Gamma_N)$.
\qed

Now we prove Theorem \ref{2gtdig-regularnormal-1} to show that: for an integer   $s\geq 2$ and  a $(G,s)$-geodesic-transitive
digraph $\Gamma $,
if $G$ has a  nontrivial regular normal subgroup, then $\Gamma$ is a circuit.

\medskip
\noindent {\bf Proof of Theorem \ref{2gtdig-regularnormal-1}.}
Let $(u,v)$ be an  arc. Then  $\Gamma^+(u)\neq \Gamma^+(v)$.
Suppose that   $ \Gamma^+(u)\cap \Gamma^+(v) =  \emptyset$. Then each 2-arc of $\Gamma$ is a 2-geodesic.
Since $\Gamma $ is  $(G,s)$-geodesic-transitive
for some  $s\geq 2$, it follows that $\Gamma $ is   $(G,2)$-arc-transitive, and  by \cite[Theorem 3.1]{Praeger-1989}, $\Gamma$ is a directed cycle.

In the remainder, we consider the case  that    $ \Gamma^+(u)\cap \Gamma^+(v) \neq  \emptyset$.
Assume   that $\Gamma$ is not a directed cycle.
Then  the valency  $m$ of  $\Gamma$ is at least 2.

%Assume that $\Gamma$ is not a directed cycle. Then  $\Gamma$ has valency $m\geq  2$.

Since $N$ is regular  on $V(\Gamma)$, we can identify $V(\Gamma)$ with $N$ so
that $\Gamma=\Cay(N,S)$ where $S$ is a subset of $N\setminus \{1_N\}$ and $|S|=m\geq 2$. Moreover,
$N$ acts by right multiplication and, for $x\in N$ and $i\geq  1_N$, we denote by $\Gamma_i^+(x)$ the set of
vertices at distance $i$ from $x$.

%Let $x\in \Gamma^+(1)$. Then as $G$ is transitive on 2-geodesics we must have that $G_1\cap G_x$ transitive on the
%set  of 2-geodesics of the form $(1, x, y)$ for some $y\in \Gamma^+(x)\setminus \Gamma^+(1)$. Translating the arc
%$(1, x)$ by $x\in N$ we see that $(x, x^2)$ is an arc, that is, $(1, x, x^2)$ is a 2-arc.

Since $N$ is a normal subgroup of $G$, it follows that $G_{1_N}=\Aut(N,S)$, and by the $(G,s)$-geodesic-transitivity with $s\geq 2$,
$\Aut(N,S)$
is transitive on both $S=\Gamma^+(1_N)$ and $\Gamma_2^+(1_N)$.
Thus all elements of $S$ have the same order,  and all elements of
$\Gamma_2^+(1_N)$ have the same order.

Let $x\in S$. Assume that $\langle x\rangle \setminus \{1_N\} \nsubseteq S$.
Let $i$ be the smallest
positive integer such that $x^i \notin  S$. Note that $x^i\neq  1_N$ and due to  $\langle x\rangle \setminus \{1_N\} \nsubseteq S$, we know that
$x^i = x \times x^{i-1}\in  \Gamma_2^+(1_N) \cap  \Gamma^+(x)$.
As $\Gamma$ is $(G,2)$-geodesic-transitive, $G_{1_N,x}$ acts
transitively on $\Gamma_2^+(1_N)\cap \Gamma^+(x)$. Since  $G_{1_N,x} \leq
G_{1_N} \leq \Aut(N,S)$ and $G_{1_N,x}$ fixes $x$, it follows that
$G_{1_N,x}$ fixes $x^i$. Hence  $\Gamma_2^+(1_N)\cap \Gamma^+(x)=\{x^i\}$,
and so $\Gamma^+(x)=(\Gamma^+(x)\cap S)\cup \{x^i\}$. It concludes that
$|\Gamma_2^+(1_N)\cap \Gamma^+(x)|=1$, and $|\Gamma^+(1_N)\cap \Gamma^+(x)|=m-1$.
Since  $G_{1_N}=\Aut(N,S)$
is transitive on $S$ and $|S|=m$, it follows that $[\Gamma^+(1_N)]$ induces an undirected subgraph, which is impossible.

Thus  $\langle x\rangle \setminus \{1_N\} \subseteq S$.
Assume that $o(x)=r\geq 3$. Then $x^{r-1}\neq x$ and $x^{r-1}\in S$. It leads to that  $(x,x\times x^{r-1}=1_N)$ is an arc, and so both
$(1_N,x)$ and $(x,1_N)$ are arcs. Since $\Gamma$ is $G$-arc-transitive, it follows  that $\Gamma$ is an undirected graph, a contradiction. As a consequence
the order of $x$ must be $2$.
Since all elements of $S$ have the same order, we know  that $o(a)=2$  for all $a\in S$, so  $S=S^{-1}$, and this fact forces that
$\Gamma$ is an undirected graph, which is a contradiction.

Therefore, $\Gamma$ has valency 1 and it is a directed cycle.

Finally, assume  $\Gamma =\Cay(T,S)$ is  a $(G,s)$-geodesic-transitive
normal Cayley digraph for some  $s\geq 2$.
Then $T$ is  a
nontrivial normal subgroup of $G$ which acts regularly  on $V(\Gamma)$.
Then by the previous argument,  $\Gamma$ is a directed cycle.
The proof is completed. \qed

Let $\Gamma $ be a connected $(G,2)$-geodesic-transitive
but not $(G,2)$-arc-transitive digraph. If  $G$ is    quasiprimitive  on $V(\Gamma_N)$ of type HA, HS, HC or  TW, then $G$ has a normal subgraph that acts regularly on the vertex set, and  by Theorem \ref{2gtdig-regularnormal-1},
$\Gamma$ is a directed cycle.

\medskip
By  Theorem  \ref{2gtdig-quotient-1},  for $s\geq 2$ each connected $(G, s)$-geodesic-transitive digraph has a connected  $(G, s)$-geodesic-transitive quotient digraph
corresponding to a normal subgroup $N$ of $G$ such that $G/N$ acts quasiprimitively or bi-quasiprimitively on the vertex set of the quotient digraph. Thus a preliminary step in determining all $s$-geodesic-transitive digraphs may be the determination of the base one.
We shall achieve this in the following proposition
in the case where $G$ is soluble.

\begin{prop}\label{2gtdig-irreducible-2}
 Let $\Gamma $ be a connected  $(G, 2)$-geodesic-transitive digraph such that
$G$ is soluble.  Suppose that $G$ acts quasiprimitively or bi-quasiprimitively on $V(\Gamma)$. Then $\Gamma $ is a circuit with $r$ vertices  where $r = 4$ or
$r$ is a prime.

\end{prop}
\proof Since $G$ is soluble, it follows that $G$ has a non-trivial abelian characteristic subgroup, and say $N$.
As $G$ acts quasiprimitively or bi-quasiprimitively on $V(\Gamma)$,
$N$ has at most two orbits in the vertex set.

If $N$ has  one orbit in $V(\Gamma) $, then  by Theorem  \ref{2gtdig-regularnormal-1}, $\Gamma$ is a circuit with $r$ vertices and
$r$ must be 4 or a prime.

In the remainder we  assume that $N$ has exactly  two orbits in $V(\Gamma) $.
Applying  Lemma \ref{2gtdig-normalsubg-bip},
 each 2-arc of $\Gamma$ is a 2-geodesic.
Since $\Gamma $ is   $(G, 2)$-geodesic-transitive,
it follows that $\Gamma $ is  a $(G, 2)$-arc-transitive digraph.

Suppose first  that   $N$ is  semiregular on $V(\Gamma)$.
Then by Theorem 3.3 of \cite{Praeger-1989}, $\Gamma$ is circuit with $r$ vertices where
$r$ must be 4 or a prime.

Assume now that $N$ is not semiregular on $V(\Gamma)$. Let $\Delta_1$ and $\Delta_2$ be the  two orbits of  $N$. Let $a\in \Delta_1$, and
$b\in \Gamma_1^+(a)\subseteq \Delta_2$.
Since $N$ is abelian, stabilizers $N_a$ and $N_b$ are  normal subgroups of  $N$.
Moreover, $ N_a$  fixes $\Delta_1$ pointwise and fixes $\Gamma_1^+(a)$ setwise. Thus all $N_a$ -orbits in $\Delta_2$ have length $k$ which is  a divisor  of the valency of $\Gamma$.

Similarly, the vertex stabilizer $ N_b$ fixes $\Delta_2$ pointwise,  and has orbits in $\Delta_1$ of length $k$.

Further, as $N$ is a  characteristic subgroup of $G$, it follows that  $N_aN_b$ is a non-trivial abelian normal subgroup of $G$ with at least two orbits in $V(\Gamma)$
and all orbits of length $k$.
Since $G$ acts bi-quasiprimitively on $V(\Gamma)$, $N_a N_b$ has exactly two orbits and this
forces $\Gamma$ to be an undirected complete bipartite graph which is a contradiction. We conclude  the proof. \qed

\section{Two-geodesic-transitive digraphs of small valency}

In this section,  we investigate the relationship of    $(G,2)$-geodesic-transitive property and $(G,2)$-arc-transitive property of  small valency digraphs.

%Throughout this section, we assume that $\Gamma$ is a $G$-arc-transitive
%digraph of valency $r\geq 2$.

We first prove the following lemma:

\begin{lemma}\label{val-common}
Let  $\Gamma$ be a $G$-arc-transitive
digraph of valency $r\geq 2$.  Then $|\Gamma^+(u)\cap \Gamma^+(v)|\neq r-1$ for each  arc $(u,v)$.

\end{lemma}
\proof Let $(u,v)$ be an arc.
Suppose that $|\Gamma^+(u)\cap \Gamma^+(v)|=r-1$.
Then since $\Gamma$ is $G$-arc-transitive, $G_u$ is transitive on
$\Gamma^+(u)$, and so $[\Gamma^+(u)]$ is a vertex-transitive digraph with valency $r-1$ and $r$ vertices.
Set $\Gamma^+(u)=\{v=v_1,v_2,v_3,\ldots,v_r\}$ and assume $v_1\rightarrow v_2$, $v_1\rightarrow v_3,\cdots$, and $v_1\rightarrow v_r$.
Due to  $|\Gamma^+(u)\cap \Gamma^+(v_2)|=r-1$, we must have $v_2\rightarrow v_1$, which is a contradiction.
Thus $|\Gamma^+(u)\cap \Gamma^+(v)|\neq r-1$.  \qed

Lemma \ref{val-common} leads  directly to the following result about valency 2 digraphs.

\begin{lemma}\label{val-2-1}
Let  $\Gamma$ be a  digraph of valency $2$.  If  $\Gamma$ is $(G,2)$-geodesic-transitive, then    $\Gamma$ is $(G,2)$-arc-transitive.

\end{lemma}
\proof   Suppose  that $\Gamma$ is $(G,2)$-geodesic-transitive.
Let $(u,v)$ be an arc. Then since $\Gamma$ has valency $2$ and applying  Lemma \ref{val-common}, we have $|\Gamma^+(u)\cap \Gamma^+(v)|\neq 1$,
it follows that  $\Gamma^+(u)\cap \Gamma^+(v)=  \emptyset$.
Thus by Lemma \ref{2gtdig-lem-at-1},  each 2-arc of $\Gamma$ is a 2-geodesic, and so $\Gamma$ is $(G,2)$-arc-transitive.
 \qed

%\subsection{Two-geodesic-transitive digraphs of valency $3$}

For  $G$-arc-transitive
digraphs of valency $3$, we have the following claim.

\begin{lemma}\label{val3-2gt-1}
Let  $\Gamma$ be a  $(G,2)$-geodesic-transitive
digraph of valency $3$.  Then  $\Gamma$ is $(G,2)$-arc-transitive.
\end{lemma}
\proof  Let $(u,v)$ be an arc of $\Gamma$.
Since  $\Gamma$ is $(G,2)$-geodesic-transitive
 of valency $3$, it follows that
$|\Gamma^+(u)\cap \Gamma^+(v)|=0,1, 2$.
Furthermore, by Lemma \ref{val-common},   $|\Gamma^+(u)\cap \Gamma^+(v)|\neq 2$.

Now consider the case  that $|\Gamma^+(u)\cap \Gamma^+(v)|=1$.
The  $G$-arc-transitive property of $\Gamma$ indicates that  $G_u$ is transitive on
$\Gamma^+(u)$, and so $[\Gamma^+(u)]$ is a vertex-transitive digraph with valency 1 and 3 vertices.
Hence  $[\Gamma^+(u)]$ is a directed circuit with 3 vertices.
 Since $\Gamma$ is $G$-vertex-transitive, for each vertex $x$,  $[\Gamma^+(x)]$ is a directed circuit with 3 vertices.

Set $\Gamma^+(u)=\{v=v_1,v_2,v_3\}$ and  $\Gamma_2^+(u)\cap \Gamma^+(v_1)=\{w_1,w_2\}$. Let $v_1\rightarrow v_2 \rightarrow v_3 \rightarrow v_1$. Then $\Gamma^+(v_1)=\{v_2,w_1,w_2\}$. For the reason that $[\Gamma^+(v_1)]$ is a directed circuit with 3 vertices, without loss of generality,
we assume that $v_2\rightarrow w_2 \rightarrow w_1 \rightarrow v_2$.

Since $\Gamma$ is  $(G,2)$-geodesic-transitive, it follows that $G_{u,v_1}$ is transitive on
$\Gamma_2^+(u)\cap \Gamma^+(v_1)=\{w_1,w_2\}$. Due to  $w_1 \rightarrow v_2$ and  $G_{u,v_1}=G_{u,v_2}=G_{u,v_3}$,
we must have  $w_2 \rightarrow v_2$, contradicts the fact that $v_2 \rightarrow w_2$.
Thus $|\Gamma^+(u)\cap \Gamma^+(v)|\neq 1$.

It concludes that   $|\Gamma^+(u)\cap \Gamma^+(v)|=0$, and $\Gamma^+(u)\cap \Gamma^+(v)=  \emptyset$. Therefore   each 2-arc of $\Gamma$ is a 2-geodesic, and so $\Gamma$ is $(G,2)$-arc-transitive. \qed

%\subsection{Two-geodesic-transitive digraphs of valency $4$}

\begin{lemma}\label{n2gt-valr-common-2}
Let  $\Gamma$ be a $(G,2)$-geodesic-transitive
digraph of valency $r\geq 3$. Let $(u,v)$ be an arc. Suppose that  $[\Gamma^+(u)]\cong k\Sigma$ where $k\geq 1$, $\Sigma$ is connected and $|V(\Sigma)|\geq 3$. Then   $|\Gamma^+(u)\cap \Gamma^+(v)|\neq 1$.

\end{lemma}
\proof
Suppose to the contrary that  $|\Gamma^+(u)\cap \Gamma^+(v)|=1$.
We assume that  $r=ke$ where  $e=|V(\Sigma)|\geq 3$. Since $\Gamma$ is $G$-arc-transitive, $G_u$ is transitive on
$\Gamma^+(u)$, and so $[\Gamma^+(u)]$ is the union of $k$ disjoint $\Sigma$ where $\Sigma$ is a connected circuit with $e$ vertices.
Moreover, the   $G$-vertex-transitive property of $\Gamma$ indicates that  for each vertex $x$,  $[\Gamma^+(x)]$ is the union of $k$ disjoint $\Sigma$.

Set $\Gamma^+(u)=\{v=v_{11},v_{12},\ldots,v_{1e},v_{21},v_{22},\ldots,v_{2e},\ldots,v_{k1},v_{k2},\ldots,v_{ke}\}$ and assume $v_{11}\rightarrow v_{12} \rightarrow \cdots   \rightarrow v_{1e} \rightarrow v_{11}$.
Then  $[v_{11},v_{12},\ldots,v_{1e}]\cong \Sigma$ is a connected circuit. Since $[\Gamma^+(u)]\cong k\Sigma$, it follows that
$G_{u,v_{11}}=G_{u,v_{12}}=G_{u,v_{13}}=\cdots=G_{u,v_{1e}}$.

As  $|\Gamma^+(u)\cap \Gamma^+(v_{11})|=1$, we have  $|\Gamma_2^+(u)\cap \Gamma^+(v_{11})|=r-1=ke-1=(k-1)e+(e-1)$.
Set  $\Gamma_2^+(u)\cap \Gamma^+(v_{11})=\{w_{11},w_{12},\ldots,w_{1(e-1)},w_{21},w_{22},\ldots,w_{2e},\ldots,w_{k1},w_{k2},\ldots,w_{ke}\}$.  Then $\Gamma^+(v_{11})=\{v_{12},w_{11},w_{12},\ldots,w_{1(e-1)},w_{21},w_{22},\ldots,w_{2e},\ldots,w_{k1},w_{k2},\ldots,w_{ke}\}$.
Note that  $[\Gamma^+(v_{11})]\cong k\Sigma$. Without loss of generality, we can assume that
$v_{12}\rightarrow w_{11} \rightarrow w_{12} \rightarrow w_{13} \rightarrow \ldots \rightarrow w_{1(e-1)} \rightarrow v_{12}$.
Since $\Gamma$ is $(G,2)$-geodesic-transitive, it follows that $G_{u,v_{11}}$ is transitive on $\Gamma_2^+(u)\cap \Gamma^+(v_{11})$.
Further, applying the facts that  $G_{u,v_{11}}=G_{u,v_{12}}$ and  $ w_{1(e-1)} \rightarrow v_{12}$, we would have  $ w_{11} \rightarrow v_{12}$,
contradicts that $v_{12}\rightarrow w_{11}$.

Thus $|\Gamma^+(u)\cap \Gamma^+(v)|\neq 1$.  \qed

Note that in Lemma \ref{n2gt-valr-common-2}, if $k=1$, then
$[\Gamma^+(u)]$ is connected; and if $k>1$, then
$[\Gamma^+(u)]$ is disconnected.

\begin{lemma}\label{n2gt-val4-common-1}
Let  $\Gamma$ be a $(G,2)$-geodesic-transitive
digraph of valency $4$.  Then for each arc  $(u,v)$,   $|\Gamma^+(u)\cap \Gamma^+(v)|\neq 1$.

\end{lemma}
\proof Let $(u,v)$ be an arc. Suppose that $|\Gamma^+(u)\cap \Gamma^+(v)|= 1$.
Then by the   $G$-arc-transitive property of $\Gamma$, the vertex stabilizer $G_u$ is transitive on
$\Gamma^+(u)$, and so $[\Gamma^+(u)]$ is a vertex-transitive digraph with valency 1 and 4 vertices.
Moreover, $[\Gamma^+(u)]$ is connected, and applying    Lemma \ref{n2gt-valr-common-2},
we have $|\Gamma^+(u)\cap \Gamma^+(v)|\neq 1$, which is a contradiction.
 \qed

 \begin{lemma}\label{val-common2}
Let  $\Gamma$ be a $G$-arc-transitive
digraph of valency $r\geq 4$.  Then $|\Gamma^+(u)\cap \Gamma^+(v)|\neq r-2$ for each   arc $(u,v)$.
\end{lemma}
\proof Suppose to the contrary that $|\Gamma^+(u)\cap \Gamma^+(v)|=r-2$ for some arc  $(u,v)$.
Set $\Gamma^+(u)=\{v=v_1,v_2,v_3,\ldots,v_r\}$. We can  assume $v_1\rightarrow v_2,\ldots$, and $v_1\rightarrow v_{r-1}$.
For the reason that  $r\geq 4$, we have  $r-1\geq 3$, and so      $v_2\nrightarrow v_1$ and $v_3\nrightarrow v_1$.

Since $\Gamma$ is  $G$-arc-transitive, it follows that $G_u$ acts  transitively on
$\Gamma^+(u)$,  and so  $|\Gamma^+(u)\cap \Gamma^+(v_2)|=r-2$ and $|\Gamma^+(u)\cap \Gamma^+(v_3)|=r-2$. Hence we must have $v_2\rightarrow v_3,\ldots$, $v_2\rightarrow v_r$,
and  $v_3\rightarrow v_2$, $v_3\rightarrow v_4,\ldots,$ and $v_3\rightarrow v_r$,
which is a contradiction.
Thus $|\Gamma^+(u)\cap \Gamma^+(v)|\neq r-2$.  \qed

\begin{lemma}\label{2gt-2at-common-1}
Let  $\Gamma$ be a $(G,2)$-geodesic-transitive
digraph of valency  $4$.  Then  $\Gamma$ is  $(G,2)$-arc-transitive.

\end{lemma}
\proof Suppose that $\Gamma$ is  $(G,2)$-geodesic-transitive. For  each  arc $(u,v)$,
since $\Gamma$ has valency 4, it follows that  $|\Gamma^+(u)\cap \Gamma^+(v)|\leq 3$, that is, $|\Gamma^+(u)\cap \Gamma^+(v)|=0,1,2,3$.
Moreover, by Lemmas \ref{val-common},  \ref{n2gt-val4-common-1} and \ref{val-common2},  $|\Gamma^+(u)\cap \Gamma^+(v)|\neq 1,2,3$.
Thus
we must have  $\Gamma^+(u)\cap \Gamma^+(v)=0$.
It leads to that each 2-arc is a 2-geodesic, and so $\Gamma$ is  $(G,2)$-arc-transitive. \qed

%\subsection{Two-geodesic-transitive digraphs of valency $5$}

\begin{lemma}\label{val5-common}
Let  $\Gamma$ be a $G$-arc-transitive
digraph of valency $5$. If    $|\Gamma^+(u)\cap \Gamma^+(v)|=2$ for some arc $(u,v)$, then $\Gamma$ is not  $(G,2)$-geodesic-transitive.
\end{lemma}
\proof Let $(u,v)$ be an arc.
Suppose that $|\Gamma^+(u)\cap \Gamma^+(v)|=2$.
Since $\Gamma$ is $G$-arc-transitive, $G_u$ is transitive on
$\Gamma^+(u)$, and so $[\Gamma^+(u)]$ is a vertex-transitive digraph with valency $2$ and $5$ vertices.
Moreover, for each vertex $x\in \Gamma^+(u)$, the set  $\Gamma^+(u)\cap \Gamma^+(x)$ contains precisely $2$ vertices.

Set $\Gamma^+(u)=\{v=v_1,v_2,v_3,v_4,v_5\}$ and assume
\begin{center}
$v_1\rightarrow v_2$, $v_1\rightarrow v_3$.
\end{center}

Suppose that $v_2\nrightarrow v_3$ and $v_3\nrightarrow v_2$.
Then $v_2\rightarrow v_4$, $v_2\rightarrow v_5$ and $v_3\rightarrow v_4$, $v_3\rightarrow v_5$.
As a consequence, $v_4\nrightarrow v_2$,  $v_4\nrightarrow v_3$ and $v_5\nrightarrow v_2$,  $v_5\nrightarrow v_3$.
Hence  $v_4\rightarrow v_1$ and $v_5\rightarrow v_1$.
Since $|\Gamma^+(u)\cap \Gamma^+(v_4)|=2$ and  $v_4\nrightarrow v_2$,  $v_4\nrightarrow v_3$,  it follows that  $v_4\rightarrow v_5$, and so
$v_5\nrightarrow v_4$.
By the previous, $v_5\nrightarrow v_2$ and   $v_5\nrightarrow v_3$, we must have $\Gamma^+(u)\cap \Gamma^+(v_5)=\{v_1\}$, which contradicts that $|\Gamma^+(u)\cap \Gamma^+(v_5)|=2$.

Thus either $v_2\rightarrow v_3$ or  $v_3\rightarrow v_2$. Without loss of generality, assume that
$$v_2\rightarrow v_3.$$
Then $v_3\nrightarrow v_2$, and as a result  $v_3\rightarrow v_4$ and $v_3\rightarrow v_5$. Furthermore,  $v_2 \rightarrow$  one of $v_4,v_5$. Assume that $v_2\rightarrow v_4$.
Then $v_4\nrightarrow v_2$ and $v_4\nrightarrow v_3$. Hence
$v_4\rightarrow v_1$ and  $v_4\rightarrow v_5$.

%Set $\Gamma_2^+(u)\cap \Gamma^+(v_1)=\{w_1,w_2,w_3\}$.
%Then $\Gamma^+(v_1)=\{v_2,v_3,w_1,w_2,w_3\}$ and $[\Gamma^+(v_1)]$ is a digraph with valency $2$ and $5$ vertices.

Since $\Gamma^+(u)\cap \Gamma^+(v_1)=\{v_2,v_3\}$ and since  $v_2\rightarrow v_3$,    $v_3\nrightarrow v_2$, it follows that
$G_{u,v_1}$ also fixes $v_2$ and $v_3$ pointwise, and so $G_{u,v_1}=G_{u,v_2}=G_{u,v_3}$.
Due to    $v_4\rightarrow v_5$, and by a similar argument, we would have  $G_{u,v_1}=G_{u,v_4}=G_{u,v_5}$.
Thus
\begin{center}
$G_{u,v_1}=G_{u,v_2}=G_{u,v_3}=G_{u,v_4}=G_{u,v_5}$. \ \ \ \ \ \ \ \ \ \ $(*)$
\end{center}

It leads to the fact   $G_{u,v_1}^{\Gamma^+(u)}=1$, that is,   $G_{u}^{\Gamma^+(u)}$ acts regularly on $\Gamma^+(u)$.
 Let $K$ be the kernel of the
$G_{u}$-action   on $\Gamma^+(u)$.
Then $G_{u}^{\Gamma^+(u)}\cong G_u/K$.

Suppose that $\Gamma$ is $(G,2)$-geodesic-transitive, then $G_{u,v_1}$ is transitive on
$\Gamma_2^+(u)\cap \Gamma^+(v_1)$.
Since $G_{u}^{\Gamma^+(u)}$ acts regularly on $\Gamma^+(u)$ and  $G_{u}^{\Gamma^+(u)}\cong G_u/K$, it follows that
$K$ is transitive on
$\Gamma_2^+(u)\cap \Gamma^+(v_1)$.
However,  each element  $k\in K$  fixes vertices $v_1,v_2$ and $v_3$, and so      $k\in G_{v_1,v_2,v_3}$.
For the reason that  $v_2,v_3\in  \Gamma^+(v_1)$,  applying the $G$-arc-transitive property of $\Gamma$ and $(*)$, we know that  $k$ acts trivially on
$\Gamma^+(v_1)$, and as a consequence
$k$ acts trivially  on
$\Gamma_2^+(u)\cap \Gamma^+(v_1)$, which is a contradiction. Therefore  $\Gamma$ is not  $(G,2)$-geodesic-transitive. \qed

\begin{lemma}\label{val5-lem2}
Every $(G,2)$-geodesic-transitive
digraph of valency $5$ is also    $(G,2)$-arc-transitive.
\end{lemma}
\proof Let  $\Gamma$ be a   $(G,2)$-geodesic-transitive digraph of valency 5 and  let $(u,v)$ be an arc.
Then   $|\Gamma^+(u)\cap \Gamma^+(v)|\leq 4$, that is, $|\Gamma^+(u)\cap \Gamma^+(v)|=0,1,2,3,4$.
Moreover, by Lemmas \ref{val-common}, \ref{val-common2} and \ref{val5-common},  $|\Gamma^+(u)\cap \Gamma^+(v)|\neq 2,3,4$.

Assume $|\Gamma^+(u)\cap \Gamma^+(v)|= 1$. Since $\Gamma$ is $G$-arc-transitive, $G_u$ is transitive on
$\Gamma^+(u)$, and so $[\Gamma^+(u)]$ is a vertex-transitive digraph with valency 1 and 5 vertices.
It follows from Lemma \ref{n2gt-valr-common-2} and the fact  $[\Gamma^+(u)]$ is connected that
$|\Gamma^+(u)\cap \Gamma^+(v)|\neq 1$, which is a contradiction.

Thus we must have  $\Gamma^+(u)\cap \Gamma^+(v)=0$.
It leads to that each 2-arc is a 2-geodesic, and so $\Gamma$ is  $(G,2)$-arc-transitive. \qed

Now we prove Theorem \ref{2gtdig-smallval-th1}.

\medskip
\noindent {\bf Proof of Theorem \ref{2gtdig-smallval-th1}.}
Let  $\Gamma$ be a $G$-arc-transitive
digraph of valency $r\geq 2$.
Since each 2-geodesic of $\Gamma$ is a 2-arc, it follows that $\Gamma$ is $(G,2)$-arc-transitive indicating that
it is also $(G,2)$-geodesic-transitive.
Assume conversely that  $\Gamma$ is  $(G,2)$-geodesic-transitive.
If $\Gamma$ has   valency  $r$ where  $2\leq r \leq 5$, then by Lemmas \ref{val-2-1}, \ref{val3-2gt-1}, \ref{2gt-2at-common-1}, \ref{val5-lem2},
 $\Gamma$ is  $(G,2)$-arc-transitive, and (i) holds.

Now let  $\Gamma$ be a $(G,2)$-geodesic-transitive digraph of diameter $2$.
Then $\Gamma$ is $G$-distance-transitive. By \cite[Theorem 4.2]{Lam-1980}, $\Gamma$ is known, and it is  a balanced incomplete block design with the Hadamard parameters, and (ii) holds.
\qed


\begin{thebibliography}{hhhh}

\bibitem{ACMX-1996}
B. Alspach, M. Conder, D. Maru$\check{{\rm s}}$i$\check{{\rm c} }$ and M. Y. Xu, A
classification of 2-arc-transitive circulants, {\it J. Algebraic
Combin.} {\bf 5} (1996), 83--86.

%\bibitem{Baddeley-1}
%R. W. Baddeley,  Two-arc transitive graphs and twisted wreath
%products, {\it J. Algebraic Combin.} {\bf 2} (1993), 215--237.

%\bibitem{BSP04-1}
%R. W. Baddeley, C. E. Praeger and C. Schneider, Transitive simple subgroups of wreath products in product action,
%{\it J. Aust. Math. Soc.} {\bf 77(1)} (2004), 55--72.

%\bibitem{Bannai-1972}
%E. Bannai,  Maximal subgroups of low rank of finite symmetric and
%alternating groups. {\it J. Fac. Sci. univ. Tokyo, Sect.} {\bf IA
%18} (1972), 475--486.


%\bibitem{BCN}
%A. E. Brouwer, A. M. Cohen and A. Neumaier, Distance-Regular Graphs,
%Springer Verlag, Berlin, Heidelberg, New York, (1989).

%\bibitem{BW-1990}
%A. E. Brouwer and  H. A. Wilbrink, Ovoids and fans in the generalized quadrangle $GQ(4, 2)$,
%{\it Geom. Dedicata} 36 (1990), 121--124.

\bibitem{Cameron-1}
P. J. Cameron,  Permutation Groups, volume 45 of London Mathematical
Society Student Texts, Cambridge University Press, Cambridge,
(1999).

\bibitem{CPW-1993}
P. J. Cameron, C. E. Praeger and N. C. Wormald, Infinite highly arc transitive digraphs and universal covering digraphs, {\it Combinatorica} {\bf 13} (1993), 377--396.



%\bibitem{Cohen-1}
%Arjeh M. Cohen,  Local recognition of graphs, buildings, and related
%geometries. In Finite Geometries, Buildings, and related Topics
%(edited by William M. Kantor, Robert A. Liebler, Stanley E. Payne,
%Ernest E. Shult),  Oxford Sci. Publ., New York. {\bf 19} (1990), 85--94.

\bibitem{CLP-1995}
M. Conder, P. Lorimer and C. Praeger, Constructions for arc-transitive digraphs, {\it J. Austral. Math. Soc. Ser. A} {\bf 59} (1995), 61--80.



%\bibitem{DDP-1}
%A. Daneshkhah, A. Devillers and C. E. Praeger, Symmetry properties
%of subdivision graphs, {\it Discrete Math.} (2011), doi:10.1016/j.disc. 2011.03.031.

%\bibitem{DGLPP-rank3-2011}
%A. Devillers, M. Giudici, C. H. Li, G. Pearce and C. E. Praeger, On
%imprimitive rank 3 permutation groups, {\it J. London Math. Soc.}.

%\bibitem{DGLP-locdt-2012}
%A. Devillers, M. Giudici, C. H. Li and C. E. Praeger, Locally
%$s$-distance transitive Graphs, {\it J. Graph Theory } {\bf 69(2)}
%(2012), 176-197.

\bibitem{DJLP-clique}
A. Devillers, W. Jin, C. H. Li and C. E. Praeger, Local $2$-geodesic
transitivity and clique graphs, {\it J. Combin. Theory Ser. A} {\bf
120} (2013), 500--508.

\bibitem{DJLP-2}
A. Devillers, W. Jin, C. H. Li and C. E. Praeger, Line graphs and
geodesic transitivity, {\it Ars Math. Contemp.} {\bf 6} (2013),
13--20.


\bibitem{DJLP-prime}
A. Devillers, W. Jin, C. H. Li and C. E. Praeger, Finite 2-geodesic
transitive graphs of prime valency, {\it J. Graph Theory} {\bf
80} (2015), 18--27.

\bibitem{DM-1}
J. D. Dixon and B. Mortimer,  Permutation groups, Springer, New
York, (1996).

\bibitem{Evans-1997}
D. M. Evans, An infinite highly arc-transitive digraphs, {\it Europ. J. Combin.} {\bf
18} (1997), 281--286.

%\bibitem{FLP-1}
%X. G. Fang, C. H. Li and C. E. Praeger, The locally two-arc
%transitive graphs admitting a Ree simple group,  {\it J. Algebra}
%{\bf 282} (2004),  638--666.

\bibitem{FH-2018}
R. Q. Feng and P. C. Hua, A new family of geodesic transitive graphs, {\it Discrete. Math.}, {\bf 341} (2018), 2700--2707.

\bibitem{GLP-1}
M. Giudici, C. H. Li and C. E. Praeger, Analysing finite locally
$s$-arc transitive graphs, {\it Trans. Amer. Math. Soc.} {\bf 356}
 (2003), 291--317.

%\bibitem{GLP-2006}
%M. Giudici, C. H. Li and C. E. Praeger, Characterising finite
%locally $s$-arc transitive  graphs with a star normal quotient, {\it
%J. Group Theory} {\bf 9(5)} (2006),  641--658.

\bibitem{GLX-2017}
M. Giudici, C. H. Li and B. Z. Xia, An infinite family of vertex-primitive
$2$-arc-transitive  digraphs, {\it
J. Combin. Theory Ser. B} {\bf 127} (2017),  1--13.

\bibitem{GX-2018}
M. Giudici  and B. Z. Xia, Vertex-quasiprimitive
$2$-arc-transitive  digraphs, {\it
Ars Math. Combin.} {\bf 14} (2018),  67--82.

%\bibitem{Guralnick-1}
%R. M. Guralnick,  Subgroups of prime power index in a simple group,
%{\it J. Algebra} {\bf 225} (1983), 304--311.

%\bibitem{DGH-1}
%D. G. Higman,  Intersection matrices for finite permutation groups,
%{\it J. Algebra} {\bf 6} (1967), 22--42.

%\bibitem{HNS-1}
%A. Hiraki, K. Nomura, H. Suzuki,  Distance regular  graphs of
%valency 6 and $a_1=1$, {\it J. Algebraic Combin. } {\bf 11} (2000),
%101--134.

\bibitem{HFZ-1}
J. J. Huang, Y. Q. Feng and J. X. Zhou,  Two-geodesic transitive graphs of order  $p^n$ with $n\leq 3$,
https://arxiv.org/abs/2207.10919v2.

\bibitem{IP-1}
A. A. Ivanov and C. E. Praeger,  On finite affine 2-arc transitive
graphs, {\it European J. Combin.} {\bf 14} (1993), 421--444.

%\bibitem{JW-val9}
%W. Jin, Finite $2$-geodesic-transitive graphs of valency 9, {\it Util. Math.}, to appear.

%\bibitem{JW-val2p}
%W. Jin, Finite $2$-geodesic-transitive graphs of valency twice a prime, {\it European J. Combin.} {\bf 49} (2015), 117--125.

%\bibitem{JW-val3p}
%W. Jin, Finite $2$-geodesic-transitive graphs of valency $3p$, {\it Ars Combin.}, to appear.

%\bibitem{Jin-3gt}
%W. Jin, Finite $3$-geodesic-transitive but not $3$-arc-transitive graphs, {\it Bull. Aust. Math. Soc.} {\bf 91} (2015), 183--190.

\bibitem{DJLP-compare}
W. Jin, A. Devillers,  C. H. Li and C. E. Praeger, On geodesic  transitive  graphs, {\it Discrete Math.} {\bf 338} (2015), 168--173.

%\bibitem{JT-local}
%W. Jin and L. Tan, Finite $s$-geodesic transitive graphs which are locally disconnected, submitted.



%\bibitem{LPM-prime-2009}
%C. H. Li, J. M. Pan and L. Ma,  Locally primitive graphs of prime
%power order, {\it J. Aust. Math. Soc.} {\bf 86} (2009), 111--122.

%\bibitem{LPL-prime-2011}
%C. H. Li, J. M. Pan and B. G. Lou,  Finite locally primitive abelian Cayley graphs, {\it Sci. China. Math.} {\bf (4) 54} (2011), 845--854.

%\bibitem{KL-rank3-1982}
%W. M. Kantor and R. A. Liebler, The rank 3 permutation
%representations of the finite classical groups, {\it Trans. Amer.
%Math. Soc.} {\bf 271} (1982), 1-71.

%\bibitem{Kovacs-2004}
%I. Kov\'acs,  Classifying arc-transitive circulants, {\it J.
%Algebraic Combin. } {\bf 20} (2004), 353--358.

%\bibitem{LM-prime-2009}
%Cai Heng Li and Li Ma,  Locally primitive graphs of prime power
%order, {\it J. Aust. Math. Soc.} {\bf 86} (2009), 111--122.

%\bibitem{KL-1}
%W. M. Kantor, R. A. Liebler, The rank 3 permutation representations
%of the finite classical groups, {\it Trans. Amer. Math. Soc.} {\bf
%271(1)} (1982), 1--71.

\bibitem{Lam-1980}
C. W.H. Lam, Distance-transitive digraphs, {\it Discrete Math.} {\bf 29} (1980), 265--274.

%\bibitem{Leemans-2009}
%D. Leemans, Locally $s$-arc-transitive graphs related to sporadic
%simple groups, {\it J. Algebra} {\bf 322(3)} (2009), 882--892.

%\bibitem{Licaiheng-2001}
%C. H. Li, The finite vertex-primitive and vertex-biprimitive $s$-transitive graphs for $s \geq 4$,
%{\it Trans. Amer. Math. Soc.} {\bf 353} (2001), 3511--3529.

%\bibitem{LCH-circulant-2005}
%C. H. Li, Permutation groups with a cyclic regular subgroup and arc
%transitive circulants, {\it J. Algebraic Combin.} {\bf 21} (2005),
%131--136.

%\bibitem{Liebeck-1}
%M. W. Liebeck,  The affine permutation groups of rank three, {\it
%Proc. London Math. Soc.} {\bf 54(3)} (1987), 477--516.


%\bibitem{LS-1986}
%M. W. Liebeck and J. Saxl, The finite primitive permutation groups
%of rank three, {\it Bull. London Math. Soc.} {\bf 18(2)} (1986), 165--172.

%\bibitem{Lorimer-1}
%P. Lorimer, Vertex transitive graphs: symmetric graphs of prime
%valency, {\it J. Graph Theory} {\bf 8} (1984), 55--68.

\bibitem{MMSZ-2002}
A. Malni\v c, D. Maru\v si\v{c}, N. Seifter and B. Zgrabli\' c, Highly arc transitive digraphs with no homomorphism onto Z, {\it Combinatorica} {\bf 22} (2002), 435--443.

\bibitem{MW-2000}
J. X. Meng and J. Z. Wang, A classification of 2-arc transitive circulant digraphs, {\it Discrete Math.} {\bf 222} (2000), 281--284.

\bibitem{MPV-2019}
L. Morgan,  P. Poto\v cnik and G. Verret, Arc-transitive digraphs of given out-valency and with blocks of given size, {\it J. Combin. Theory Ser. B} {\bf 137} (2019), 118--125.

\bibitem{Neumann-1977}
P. M. Neumann, Finite permutation groups, edge-coloured graphs and matrices, in {\it Topics in Group Theory
and Computation}, Academic Press, London, 1977,  82--118.

%\bibitem{Paley-1}
%R. E. A. C. Paley,  On orthogonal matrices, {\it J. Math. Phys.}
%{\bf 12} (1933), 311--320.

\bibitem{Praeger-1989}
C. E. Praeger,  Highly arc  transitive digraphs,
{\it European J. Combin.} {\bf 10} (1989), 281--292.


\bibitem{Praeger-4}
C. E. Praeger,  An O'Nan Scott theorem for finite quasiprimitive
permutation groups and an application to 2-arc transitive graphs,
{\it J. London Math. Soc.} {\bf 47(2)} (1993), 227--239.


%\bibitem{Praeger-2}
%C. E. Praeger,  Finite transitive permutation groups and finite
%vertex-transitive graphs, in {\it Graph Symmetry: Algebraic Methods
%and Applications, NATO ASI Ser. C} {\bf 497} (1997), 277--318.

\bibitem{Praeger-2003-biq}
C. E. Praeger, Finite transitive permutation groups and  bipartite vertex-transitive graphs, {\it Illinois J. Mathematics} {\bf 47} (2003) 461--475.

%\bibitem{Praeger-3}
%C. E. Praeger, J. Saxl and K. Yokohama, Distance transitive graphs
%and finite simple groups, {\it Proc. London Math. Soc.} {\bf (3) 55} (1987), 1--21.



%\bibitem{Rooij-Wilf-1}
%A. C. M. van Rooij  and H. S. Wilf,   The interchange  graph of a
%finite graph, {\it Acta Mathematica Hungarica} {\bf (3-4) 16}
%(1965), 263--269.

%\bibitem{Stellmacher}
%B. Stellmacher,  Locally $s$-transitive  graphs, {\it unpublished}.


\bibitem{Tutte-1}
W. T. Tutte,  A family of  cubical  graphs, {\it Proc. Cambridge
Philos. Soc.} {\bf 43} (1947), 459--474.


\bibitem{Tutte-2}
W. T. Tutte, On the symmetry of cubic graphs, {\it Canad. J. Math.}
{\bf 11} (1959), 621--624.


\bibitem{WS-2003}
K. S. Wang and H. Suzuki, Weakly distance-regular digraphs, {\it Discrete Math.}
{\bf 264} (2003), 225--236.

\bibitem{Weiss-1}
R. Weiss,  s-transitive graphs, {\it Colloquia Mathematica
Societatis Janos Bolyai, Algebraic methods in graph theory, szeged
(Hungary)} {\bf 25} (1978), 827--847.



\bibitem{weiss}
R. Weiss,  The non-existence of 8-transitive graphs, {\it
Combinatorica} {\bf 1} (1981), 309--311.

\bibitem{Wielandt-book}
H. Wielandt, Finite Permutation Groups, New York: Academic Press
(1964).













\end{thebibliography}
\end{document}